\documentclass[12pt]{article}

\usepackage{graphicx}
\usepackage[pdftex]{color}
\usepackage{amsmath}
\usepackage{amssymb}
\usepackage{wasysym}

\newtheorem{theorem}{Theorem}[section]
\newtheorem{corollary}[theorem]{Corollary}
\newtheorem{proposition}[theorem]{Proposition}
\newtheorem{lemma}[theorem]{Lemma}

\newtheorem{problem}[theorem]{Problem}

\newcommand{\proof}{{\em Proof.} }

\newcommand{\qed}{ {\APLbox}}

\newcommand{\Range}{\mathop{\rm Range}}

\newcommand{\tr}{\mathop{\rm tr}}
\newcommand{\Rgamma}{\Range \Gamma}

\newcommand{\JpsiKL}{{J_\Psi^{KL}}}
\newcommand{\JpsiH}{{J_\Psi^{H}}}

\newcommand{\LHG}{{\cal L}^H_\Gamma}
\newcommand{\LH}{{\cal L}^H}
\newcommand{\LKLG}{{\cal L}^{KL}_\Gamma}
\newcommand{\LKL}{{\cal L}^{KL}}

\newcommand{\R}{\mathbb{R}}
\newcommand{\C}{\mathbb{C}}
\newcommand{\D}{\mathbb{D}}
\newcommand{\T}{\mathbb{T}}

\newcommand{\SpecSpace}{{{\cal S}_+^{m\times m}(\mathbb{T})}}

\newcommand{\dth}{{\rm d}\vartheta}
\newcommand{\ejth}{{\rm e}^{{\rm j}\vartheta}}
\newcommand{\emjth}{{\rm e}^{-{\rm j}\vartheta}}
\newcommand{\suchthat}{\mbox{{\rm such that}}}
\newcommand{\Unitcircle}{{\mathbb T}}
\newcommand{\Hermitian}{\mathbb{H}}

\newcommand{\SMM} {{{\cal S}^{m\times m}_+}}
\newcommand{\SuuT} {{{\cal S}^{1\times 1}_+}(\T)}
\newcommand{\SMMT}{{\SMM(\T)}}
\newcommand{\Continuous}{{\mathcal C}(\Unitcircle;\Hermitian(m))}

\newcommand{\optimum}{{\rm o}}

\newcommand{\PhioptKL}{{\Phi_{\optimum}^{KL}}}

\newcommand{\LambdaoptKL}{{\Lambda_{\optimum}^{KL}}}

\newcommand{\PhioptH} {{\Phi_{\optimum}^{H}}}
\newcommand{\WoptH} {{W_{\optimum}^{H}}}

\newcommand{\LambdaoptH}{{\Lambda_{\optimum}^{H}}}

\newcommand{\SigmaSet}{{P_\Gamma}}
\newcommand{\SigmaSetExplicit}{{\{\Sigma\in\Rgamma\ |\ \Sigma>0\}}}

\begin{document}
\title{On the well-posedness of multivariate spectrum approximation  and convergence of high-resolution spectral estimators.
\thanks{Partially supported by the Ministry of
Education, University, and Research of Italy (MIUR), under Project
2006094843:  {\em New techniques and applications of identification and adaptive control}}}
\author{Federico Ramponi\thanks{
Institut f\"ur Automatik, ETH Z\"urich,
Physikstrasse 3, 8092 Z\"urich, Switzerland.
e-mail: {\tt ramponif@control.ee.ethz.ch}} 
\and Augusto Ferrante\thanks{
Dipartimento di Ingeneria dell'Informazione,
Universit\`{a} di Padova,
via Gradenigo 6/B,
I-35131 Padova,
Italy.
e-mail: {\tt augusto@dei.unipd.it}}
\and
Michele Pavon\thanks{Dipartimento di Matematica Pura e Applicata, Universit\`a di Padova, via Trieste 63, 35131 Padova, Italy.
e-mail: {\tt pavon@math.unipd.it}}
}
\maketitle
\begin{abstract}
In this paper, we establish the well-posedness of the generalized moment problems recently studied by Byrnes-Georgiou-Lindquist and coworkers, and by Ferrante-Pavon-Ramponi. We then apply these continuity results to prove almost sure convergence of a sequence of high-resolution spectral estimators indexed by the sample size.
\end{abstract}

\newpage

\section{Introduction}

Consider a linear, time invariant system
\begin{equation}
\label{SYSTEM}
x(t+1) = Ax(t) + By(t),\qquad A\in \C^{n\times n}, B\in\C^{n\times m},
\end{equation}
with transfer function
\begin{equation}
\label{formulaG}
G(z)=(zI-A)^{-1}B, 
\end{equation}
where $A$ is a stability matrix, $B$ is full column rank, 
and $(A,B)$ is a reachable pair.  
Suppose that the system is fed with
a $m$-dimensional, 
zero-mean,
wide-sense stationary process $y$ having spectrum 
$\Phi$.
The asymptotic state covariance $\Sigma$ of the system (\ref{SYSTEM}) satisfies:
\begin{equation}
\label{MOMENT_PROBLEM}
\Sigma = \int G\Phi G^\ast.
\end{equation}
Here and in the following, $G^\ast(z) = G^\top(z^{-1})$,
and integration takes place over the unit circle
with respect to the normalized Lebesgue measure $\dth/2\pi$.
%
%
Let $\SMMT$ be the family of bounded, coercive, 
$\C^{m\times m}$-valued spectral density functions on the unit
circle. Hence, $\Phi\in \SMMT$ if and only if $\Phi^{-1}\in \SMMT$.
Given a Hermitian and positive-definite $n\times n$ matrix $\Sigma$,
consider the problem of finding $\Phi \in \SMMT$ 
that satisfies (\ref{MOMENT_PROBLEM}),
i.e., that is compatible with $\Sigma$. This is a particular case of a {\em moment problem.}
In the last ten years, much research has been produced, mainly by the Byrnes-Georgiou-Lindquist school, on generalized moment problems \cite{BYRNES_GUSEV_LINDQUIST_CONVEX},
\cite{GEORGIOU_MAXIMUMENTROPY},
\cite{BYRNES_LINDQUIST_THEGENERALIZED},
\cite{GEORGIOU_RELATIVEENTROPY},
\cite{GEORGIOU_LINDQUIST_KULLBACKLEIBLER},
and analytic interpolation with complexity constraint \cite{BLOMQVIST_LINDQUIST_NAGAMUNE_MATRIXVALUED}, and their applications to spectral estimation \cite{BYRNES_GEORGIOU_LINDQUIST_THREE}, \cite{GEORGIOU_LINDQUIST_CONVEXARMA}, \cite{RAMPONI_FERRANTE_PAVON_GLOBALLYCONVERGENT} and robust control \cite{GEORGIOU_LINDQUIST_REMARKSDESIGN}.
It is worth recalling that two fundamental problems of control theory,
namely the {\em covariance extension problem} and the
{\em Nevanlinna-Pick interpolation problem of robust control}, can be recast in this form
\cite{GEORGIOU_LINDQUIST_KULLBACKLEIBLER}.

Equation (\ref{MOMENT_PROBLEM}), where the unknown is $\Phi$,
is also
a typical example of an {\em inverse} problem.
%
%
%
%
%
Recall that a problem is said to be {\em well posed}, in the sense of Hadamard,
if it admits a solution, such a solution is unique, and the solution depends 
continuously on the data.
%
%
Inverse problems are typically {\em not} well posed.
In our case, there may well be no solution $\Phi$, and when a solution exists,
there may be (infinitely) many.
It was shown in
\cite{GEORGIOU_STATECOV}, that the set of solutions is nonempty if and only if there exists
$H \in \C^{m\times n}$ such that
\begin{equation}
\label{GEORGIOU-CONDITION}
\Sigma - A\Sigma A^\ast = BH + H^\ast B^\ast.
\end{equation}
%

\noindent
%
When (\ref{GEORGIOU-CONDITION}) is feasible with $\Sigma>0$,
there are infinitely many solutions $\Phi$ to (\ref{MOMENT_PROBLEM}).
%
To select a particular solution it is natural to introduce an
optimality criterion. For control applications, however, it is  desirable that such a solution
be of limited complexity.
It should namely be rational and with an {\em a priori} bound
on its MacMillan degree.
%
One of the great accomplishments
of the Byrnes-Georgiou-Lindquist approach is having shown that
the minimization of certain {\em entropy-like functionals} leads to solutions 
that satisfy this requirement.
In \cite{GEORGIOU_STATECOV}, Georgiou
provided an explicit expression for the spectrum $\hat{\Phi}$
that exhibits {\em maximum entropy rate} among the solutions of 
(\ref{MOMENT_PROBLEM}).

\noindent

Suppose now that some {\em a priori} information about $\Phi$ is available in the form
of a spectrum $\Psi \in \SMMT$.
Given $G$, $\Sigma$, and $\Psi$, we now seek a spectrum $\Phi$,
which is closest to $\Psi$ in a certain metric,
among the solutions of (\ref{MOMENT_PROBLEM}).
Paper \cite{GEORGIOU_LINDQUIST_KULLBACKLEIBLER}
deals with such an optimization problem in the case when $y$ is a scalar process.
The criterion there is the Kullback-Leibler pseudo-distance from $\Psi$
to $\Phi$. A drawback of this approach is that it does not seem to generalize to the multivariable case.
%
This motivated us to
provide a suitable extension of the so-called {\em Hellinger distance}
with respect to which the {\em multivariable} version of the problem
is solvable (see
\cite{FERRANTE_PAVON_RAMPONI_HELLINGERVSKULLBACK} and 
\cite{RAMPONI_FERRANTE_PAVON_GLOBALLYCONVERGENT}).


The main result of this paper is contained in
Section \ref{SECTION_WELL_POSEDNESS}. We show there that, under the feasibility assumption, 
the solution to the spectrum approximation problem 
with respect to both the scalar Kullback-Leibler pseudo-distance
and the multivariable Hellinger distance
depends continuously on $\Sigma$, thereby proving that these problems are well-posed.
In Section \ref{SECTION_CONSISTENCY}
we deal with the case when only an {\em estimate}
$\hat{\Sigma}$ of $\Sigma$ is available. 
By applying the continuity results of Section \ref{SECTION_WELL_POSEDNESS},
we prove a consistency result for the  solutions
to both approximation problems.

\section{Spectrum approximation problems}
\label{SECTION_SPECTRUM_APPROX}

In this section, we 
collect some background material on spectrum approximation problems.
The reader is referred
to \cite{GEORGIOU_STATECOV}, \cite{GEORGIOU_LINDQUIST_KULLBACKLEIBLER},
\cite{FERRANTE_PAVON_RAMPONI_HELLINGERVSKULLBACK} and 
\cite{RAMPONI_FERRANTE_PAVON_GLOBALLYCONVERGENT}
for a more detailed treatment.

\subsection{Feasibility of the moment problem}

Let $\Hermitian(n)$ be the space of Hermitian $n\times n$ matrices, and
$\Continuous$ the
space of $\Hermitian(m)$-valued continuous functions defined on the unit circle.
Let the operator $\Gamma: \Continuous \rightarrow \Hermitian(n)$
be defined as follows:
\begin{equation}
\label{GAMMA}
\Gamma(\Phi) := \int G \Phi G^\ast.
\end{equation}
Consider now the {\em range} of the operator $\Gamma$
(as a vector space over the reals).
We have the following result (see \cite{RAMPONI_FERRANTE_PAVON_GLOBALLYCONVERGENT}).
\begin{proposition} 
\label{RANGEGAMMAPROP} 
\ 
\begin{enumerate}
\item Let $\Sigma = \Sigma^\ast > 0$. The following are equivalent:
\begin{itemize}
\item There exists $H \in \C^{m \times n}$ which solves (\ref{GEORGIOU-CONDITION}).
\item There exists $\Phi \in \SpecSpace$ such that $\int G \Phi G^\ast = \Sigma$.
\item There exists $\Phi \in \Continuous$, $\Phi > 0$ such that $\Gamma(\Phi) = \Sigma$.
\end{itemize}
\item Let $\Sigma = \Sigma^\ast$ (not necessarily definite).
There exists $H \in \C^{m \times n}$ that solves (\ref{GEORGIOU-CONDITION})
if and only if $\Sigma \in \Rgamma$.
\item $X \in \Rgamma^\perp$ if and only if $G^\ast(\ejth)X G(\ejth) = 0 \ \forall \vartheta \in [0, 2\pi]$.
\end{enumerate}
\end{proposition}
%
%
%
We define
\begin{equation}
\SigmaSet := \SigmaSetExplicit.
\end{equation}
In view of Proposition \ref{RANGEGAMMAPROP}, 
for each $\Sigma\in\SigmaSet$
problem (\ref{MOMENT_PROBLEM}) is feasible.

\subsection{Scalar approximation in the Kullback-Leibler pseudo-distance}

In \cite{GEORGIOU_LINDQUIST_KULLBACKLEIBLER},
the Kullback-Leibler pseudo-distance for spectral densities in $\SuuT$ was introduced:
\begin{equation}
\label{KULLBACKLEIBLER}
\D(\Psi\|\Phi)=\int \Psi\log \frac{\Psi}{\Phi}.
\end{equation}
As is well known, the corresponding quantity for probability
densities originates in hypothesis testing, where it represents 
the mean information per observation for discrimination of an underlying probability density
from another \cite{KULLBACK_INFORMATIONSTATISTICS}.
The approximation problem goes as follows:
\begin{problem}
\label{PROBLEM_MINIMIZE_KL}
Given $\Sigma\in\SigmaSet$ and $\Psi\in \SuuT$, find $\PhioptKL$ that solves
\begin{equation}
\label{PROBLEM_KL_1}
\begin{split}
&{\rm minimize} \quad \D(\Psi\|\Phi)\\
&{\rm over}     \quad\left\{\Phi\in \SuuT\ |\ \int G\Phi  G^\ast= \Sigma\right\}.
\end{split}
\end{equation} 
\end{problem}
%
Note that, following \cite{GEORGIOU_LINDQUIST_KULLBACKLEIBLER},
and differently from optimization problems that are usual in the probability setting,
we minimize (\ref{KULLBACKLEIBLER}) with respect to the {\em second} argument.
The remarkable advantage of this approach is that,
differently from optimization with respect to the first argument,
it will yield a {\em rational} solution whenever $\Psi$ is rational.
Let 
\begin{displaymath}
\LKL := \{\Lambda \in \Hermitian(n) \ |\ G^\ast \Lambda G > 0, \forall e^{i\vartheta}\in\T\}.
\end{displaymath}
For a given $\Lambda\in \LKL$, consider the {\em Lagrangian functional}
\begin{equation}
\label{KL_LAGRANGIAN}
L(\Phi;\Lambda) = \D(\Psi\|\Phi) + \left< \Lambda, \int G\Phi G^\ast - \Sigma \right>,
\end{equation}
where $\left<A,B\right> := \tr AB$ denotes the scalar product between the Hermitian
matrices $A$ and $B$.
Observe that
the term $\int G\Phi G^\ast$ between brackets belongs to $\Rgamma$ by definition,
while $\Sigma$ belongs to $\Rgamma$ by the feasibility assumption.
Hence, it is natural to restrict $\Lambda$ to $\Rgamma$, or,
which is the same, to 
\begin{displaymath}
\LKLG := \LKL \cap \Rgamma.
\end{displaymath}

%
%
%
\noindent
The functional (\ref{KL_LAGRANGIAN}) is {\em strictly convex} on $\SuuT$.
Hence, its  
{\em unconstrained} minimization with respect to $\Phi$ can be pursued
imposing that its derivative in an arbitrary direction $\delta\Phi$ is zero.
This yields the form for the optimal spectrum:
\begin{equation}
\label{KL4aa}
\PhioptKL = \frac{\Psi}{G^\ast\Lambda G}.
\end{equation}
As noted previously, inasmuch as $\Psi$ is rational $\PhioptKL$
is also rational, and with MacMillan degree less than or equal to
$2n +\deg\Psi$. Now if $\Lambda \in \LKLG$ is such that
\begin{equation}
\label{KL9}
\int G\ \frac{\Psi}{G^\ast \Lambda G}\ G^\ast=\Sigma,
\end{equation} 
that is, if $\Lambda$ is such that the corresponding optimal spectrum
$\PhioptKL$ satisfies the constraint, then (\ref{KL4aa})
is the unique solution to the constrained 
approximation problem (\ref{PROBLEM_MINIMIZE_KL}).
Finding such $\Lambda$ is the objective of the {\em the dual problem},
which is readily seen 
\cite{GEORGIOU_LINDQUIST_KULLBACKLEIBLER}
to be equivalent to
\begin{equation}
\label{DUAL_PROBLEM_KL}
{\rm minimize}\quad\{\JpsiKL(\Lambda) \ | \ \Lambda\in \LKLG \}
\end{equation}
where
\begin{equation}
\label{JPSI_KL}
\JpsiKL(\Lambda)
= - \int\Psi\log G^\ast\Lambda G + \tr \Lambda \Sigma.
\end{equation}
This is also a convex optimization problem.
{\em Existence} of a minimum is a highly nontrivial issue.
Such existence was proved in 
\cite{GEORGIOU_LINDQUIST_KULLBACKLEIBLER}
resorting to a profound topological result, and in
\cite{FERRANTE_PAVON_RAMPONI_FURTHERRESULTS}
by a less abstract argument.
%
\begin{theorem}
\label{KL_EXISTENCE_THEOREM}
The strictly convex functional $\JpsiKL$
has a unique minimum point in $\LKLG$.
\end{theorem}
The minimum point
of Theorem \ref{KL_EXISTENCE_THEOREM}
provides the optimal solution to the primal problem
\ref{PROBLEM_MINIMIZE_KL} via (\ref{KL4aa}). Differently from the primal problem,
whose domain $\SuuT$ is infinite-dimensional, the 
dual problem is finite-dimensional, hence the minimization of
$\JpsiKL$ can be accomplished with iterative numerical methods.
The numerical minimization of $\JpsiKL$ is not, however, a simple problem,
because both the functional and its gradient are unbounded
on $\LKLG$ (which is unbounded itself). Moreover, reparametrization 
of $\LKLG$ may lead to loss of convexity
(see \cite{GEORGIOU_LINDQUIST_KULLBACKLEIBLER}
and references therein).
An alternative approach to this problem
was proposed in \cite{PAVON_FERRANTE_ONGEORGIOULINDQUIST}.

\subsection{Multivariable approximation in the Hellinger distance}

In \cite{FERRANTE_PAVON_RAMPONI_HELLINGERVSKULLBACK} the {\em Hellinger
distance} between two spectral densisties $\Phi, \Psi \in \SuuT$ was introduced: 
\begin{equation}
\label{HELLINGER_SCALARE}
d_H(\Phi,\Psi):=\left[\int \left(\sqrt{\Phi}-\sqrt{\Psi}\right)^2 \right]^{1/2}.
\end{equation}
As it happens for the Kullback-Leibler case, its counterpart 
for probability 
densities is well-known in mathematical statistics. Differently from
the Kullback-Leibler case, this is a {\em bona fide} distance
(note that (\ref{HELLINGER_SCALARE}) is nothing more that the $L^2$ distance
between the square roots of $\Phi$ and $\Psi$,
and that the square roots are particular instances 
of {\em spectral factors}).
A variational analysis similar to the one we have just seen is possible
and leads to similar results. 
Let us focus directly
on the multivariable extension of (\ref{HELLINGER_SCALARE}) that was developed 
in \cite{FERRANTE_PAVON_RAMPONI_HELLINGERVSKULLBACK}. 
Given $\Phi, \Psi \in \SMMT$, we define the following quantity:
\begin{equation}
\label{HELLINGER}
\begin{split}
d_H(\Phi,\Psi)&:=
\inf\left\{\|W_\Psi-W_\Phi\|_2:\ W_\Psi, W_\Phi \in L^{m\times m}_2, \right. \\
&\left. W_\Psi W_\Psi^*=\Psi,\  W_\Phi W_\Phi^*=\Phi\right\}.
\end{split}
\end{equation} 
Observe that $d_H(\Phi,\Psi)$ is simply the $L^2$ distance
between the sets of {\em all the square spectral factors} of 
$\Phi$ and $\Psi$ respectively.
We have the following result (see \cite{FERRANTE_PAVON_RAMPONI_HELLINGERVSKULLBACK}).
\begin{theorem}
\label{HELLINGER_THEOREM}
The following facts hold true:
\begin{enumerate}
\item $d_H$ is a {\em bona fide} distance function.
\item $d_H(\Phi, \Psi)$ coincides with (\ref{HELLINGER_SCALARE}) when $\Phi$ and $\Psi$ are scalar.
\item The infimum in (\ref{HELLINGER}) is indeed a minimum.
\item For any square spectral factor $\bar{W}_\Psi$ of $\Psi$, we have:
\begin{displaymath}
d_H(\Phi,\Psi) =\inf_{W_\Phi}   \left\{ \|\bar{W}_\Psi-W_\Phi\|_2:\ W_\Phi \in L^{m\times m}_2, W_\Phi W_\Phi^*=\Phi\right \}.
\end{displaymath}
\end{enumerate}
\end{theorem}
Fact 4 says that, if we fix a spectral factor of one spectrum and minimize only among
spectral factors of the other, the result is the same.
Given $\Psi \in \SMMT$ (and $G(z)$ $n\times m$), we pose a minimization problem similar
to Problem \ref{PROBLEM_MINIMIZE_KL}:
\begin{problem}
\label{PROBLEM_MINIMIZE_H}
Given $\Sigma\in\SigmaSet$ and $\Psi\in \SMMT$, find $\PhioptH$ that solves
\begin{equation}
\label{PROBLEM_H_1}
\begin{split}
&{\rm minimize} \quad d_H(\Phi, \Psi)\\
&{\rm over}     \quad\left\{\Phi \in \SMMT\ |\ \int G\Phi  G^\ast= \Sigma \right\}.
\end{split}
\end{equation} 
\end{problem}
In view of facts 3 and 4 in Theorem \ref{HELLINGER_THEOREM}, 
once a spectral factor of $\Psi$ is fixed, the same 
problem \ref{PROBLEM_MINIMIZE_H} can be reformulated in terms of
a minimization with respect to {\em spectral factors of $\Phi$}:
%
\\

\noindent
{\em Given $\Sigma\in\SigmaSet$ and 
a spectral factor $W_\Psi$ of $\Psi \in \SMMT$, find $W_\Phi$ that solves}
\begin{equation}
\label{PROBLEM_H_2}
\begin{split}
&{\rm minimize} \quad \tr \int \left( W_\Phi - W_\Psi \right)\left( W_\Phi - W_\Psi \right)^\ast \\
&{\rm over}     \quad\left\{W_\Phi \in L_2^{m\times m}\ |\ \int GW_\Phi W_\Phi^\ast  G^\ast= \Sigma \right\}.
\end{split}
\end{equation} 
%
%
Consider the Lagrangian functional
\begin{equation}
\label{H_LAGRANGIAN}
\begin{split}
H(W_\Phi,\Lambda)
= \tr \int \left( W_\Phi - W_\Psi \right)\left( W_\Phi - W_\Psi \right)^\ast 
+ \left< \Lambda, \int G W_\Phi W_\Phi^\ast G^\ast - \Sigma \right>.
\end{split}
\end{equation}
For the same reason as before, we restrict the matrix
$\Lambda$ to $\Rgamma$. 
The functional (\ref{H_LAGRANGIAN}) is {\em strictly convex}, and its  
{\em unconstrained} minimization of (\ref{H_LAGRANGIAN}) with respect to $W_\Phi$ 
yields the following condition for the optimal spectral factor
$\WoptH$
(see \cite{FERRANTE_PAVON_RAMPONI_HELLINGERVSKULLBACK} for details):
\begin{equation}
\label{H_OPTIMAL_COND}
\WoptH - W_\Psi + G^\ast\Lambda G \WoptH = 0.
\end{equation}
In order to ensure that the corresponding spectrum is integrable
over the unit circle, we now require {\em a posteriori} that $\Lambda$
belongs to the set
\begin{displaymath}
\label{LH}
\LH = \left\{ \Lambda \in \Hermitian(n)\ | \ I + G^\ast \Lambda G > 0 \ \forall \ejth \in \Unitcircle \right\}
\end{displaymath}
or, which is the same, that it
belongs to the set
\begin{equation}
\label{LHG_SET}
\LHG := \LH \cap \Rgamma.
\end{equation}
Such restriction yields the following optimal spectral factor
and spectrum:
\begin{equation}
\label{H_PRIMALMINIMUM}
\begin{split}
\WoptH    &= (I + G^\ast \Lambda G)^{-1} W_\Psi, \\
\PhioptH  &= \WoptH \WoptH^\ast =  (I + G^\ast \Lambda G)^{-1} \Psi (I + G^\ast \Lambda G)^{-1}.
\end{split}
\end{equation}
Now if $\Lambda$ is such that
\begin{equation}
\label{H9}
\int G\ (I + G^\ast \Lambda G)^{-1} \Psi (I + G^\ast \Lambda G)^{-1} \ G^\ast=\Sigma,
\end{equation} 
then $\PhioptH$ in (\ref{H_PRIMALMINIMUM}) is the unique solution to the constrained 
approximation problem (\ref{PROBLEM_MINIMIZE_H}).
In order to find such $\Lambda$, one must solve the {\em dual problem},
which can be shown to be equivalent to
\begin{equation}
\label{DUAL_PROBLEM_H}
{\rm minimize}\quad\{\JpsiH(\Lambda)\ |\ \Lambda \in \LHG \}
\end{equation}
where
\begin{equation}
\label{JPSI_H}
\JpsiH(\Lambda)
= \tr \int (I + G^\ast \Lambda G)^{-1} \Psi + \tr \Lambda\Sigma.
\end{equation}
{\em Existence} of a minimum is again a highly nontrivial issue.
We have the following result
(see \cite{FERRANTE_PAVON_RAMPONI_HELLINGERVSKULLBACK}).
\begin{theorem}
\label{H_EXISTENCE_THEOREM}
The strictly convex functional $\JpsiH$ 
has a unique minimum point in $\LHG$. 
\end{theorem}
The minimum point  of Theorem \ref{H_EXISTENCE_THEOREM}
provides the optimal solution to the primal problem
\ref{PROBLEM_MINIMIZE_H} via (\ref{H_PRIMALMINIMUM}). It can be found
by means of iterative numerical algorithms.
The numerical minimization of $\JpsiH$ is a highly nontrivial problem,
for reasons similar to the ones concerning $\JpsiKL$.
In \cite{RAMPONI_FERRANTE_PAVON_GLOBALLYCONVERGENT}, we
propose a matricial version of the Newton algorithm that
avoids any reparametrization of $\LHG$, and proved
its global convergence.

\section{Well-posedness of the approximation problems}
\label{SECTION_WELL_POSEDNESS}

\noindent
In this section, we show that both the dual problems 
(\ref{DUAL_PROBLEM_KL}) and (\ref{DUAL_PROBLEM_H})
are well-posed, since their unique
solution is continuous with respect to a small perturbation of 
$\Sigma$.
The well-posedness of the respective
{\em primal} problem then easily follows.
All these continuity properties rely on the following basic result.
\begin{theorem}
\label{CONVEXITYTHEOREM}
Let $A$ be an open and convex subset of a finite-dimensional euclidean space $V$.
Let $f:A\rightarrow \R$ be a strictly convex function, and suppose that
a minimum point $\bar{x}$ of $f$ exists.
Then,
for all $\varepsilon >0$, 
there exists $\delta>0$ such that,
for each $p \in \R^n$, $||p||<\delta$,
the function $f_p:A\rightarrow \R$
defined as
\begin{displaymath}
f_p(x) := f(x) - \left<p, x\right> 
\end{displaymath}
admits an unique minimum point $\bar{x}_p$, and moreover
\begin{displaymath}
||\bar{x}_p - \bar{x}|| < \varepsilon.
\end{displaymath}
(Note: $f^\ast(p) := -f_p(\bar{x}_p)$ is the {\em Fenchel dual} of $f$ at $p$.)
\end{theorem}

\noindent
\proof
First, note that the minimum point $\bar{x}$ is unique,
since $f$ is strictly convex.
Let $\varepsilon > 0$, and
let 
$S(\bar{x}, \varepsilon) = \{\bar{x} + y\ |\  ||y||=\varepsilon\}$ 
denote the sphere of radius $\varepsilon$ centered in $\bar{x}$.
Let moreover 
$B(\bar{x}, \varepsilon) = \{\bar{x} + y\ |\  ||y||<\varepsilon\}$ 
denote the open ball of radius $\varepsilon$ centered in $\bar{x}$ and 
$\bar{B}(\bar{x}, \varepsilon) = \{\bar{x} + y\ |\  ||y||\leq\varepsilon\}$
its closure.
Then $\bar{B}(\bar{x}, \varepsilon) = B(\bar{x}, \varepsilon) \cup S(\bar{x}, \varepsilon)$,
$\bar{B}(\bar{x}, \varepsilon)$ and $S(\bar{x}, \varepsilon)$ are compact,
and
$S(\bar{x}, \varepsilon)$ is the boundary of $B(\bar{x}, \varepsilon)$.
Since $f$ is continuous, it admits a minimum point $\bar{x} + y_\varepsilon$
over $S(\bar{x}, \varepsilon)$. 
Since  $\bar{x}$ is the unique global minimum point of $f$, we must have $m_\varepsilon:=f(\bar{x} + y_\varepsilon)-f(\bar{x})>0$.
Then, for $||y|| = \varepsilon$ we have
\begin{equation}
\label{FPMIN1}
f(\bar{x} + y) - f(\bar{x}) \geq m_\varepsilon.
\end{equation}

\noindent
Let now $0<\delta<m_\varepsilon/\varepsilon$. For $||p||<\delta$ and $||y|| = \varepsilon$
we have
\begin{equation}
\label{FPMIN2}
\left<p, y\right> \leq\ ||p||\ ||y|| < \delta \varepsilon < m_\varepsilon
\end{equation}
where the first inequality stems from the Cauchy-Schwartz inequality.
From (\ref{FPMIN1}) and (\ref{FPMIN2}), we get for $||y|| = \varepsilon$
\begin{displaymath}
\begin{split}
f(\bar{x} + y) - f(\bar{x})\ &> \left<p, y\right> = \left<p, \bar{x} + y\right> - \left<p, \bar{x}\right> \\
f_p(\bar{x} + y)\ &>\ f_p(\bar{x})
\end{split}
\end{displaymath}
that is, 
\begin{displaymath}
f_p(x)\ >\ f_p(\bar{x})
\end{displaymath}
for each $x \in S(\bar{x}, \varepsilon)$.

\noindent
Now, since $f$ is strictly convex and hence continuous, $f_p$ is also strictly convex and continuous,
and admits a minimum point $\bar{x}_p$ over the compact set 
$\bar{B}(\bar{x}, \varepsilon)$. But it follows from the previous
considerations that such minimum cannot belong to $S(\bar{x}, \varepsilon)$.
Hence, it must belong to the open ball $B(\bar{x}, \varepsilon)$.
{\em As such, $\bar{x}_p$ is also a local minimum of $f_p$ over $A$,
but since $f_p$ is strictly convex, it is also the unique global minimum point}.
Summing up, for fixed $\varepsilon>0$, there exists $\delta>0$ such that,
if $||p||<\delta$, then $f_p$ admits an unique minimum $\bar{x}_p$ over $A$.
It follows from the previous analysis that, for sufficiently small 
$\delta$, $\bar{x}_p$ belongs to $B(\bar{x}, \varepsilon)$.
This proves the theorem.
\qed

\subsection{Well-posedness of Kullback-Leibler approximation}

Consider the dual functional (\ref{JPSI_KL}), and let us make
its dependence upon $\Sigma$ explicit:
\begin{displaymath}
\JpsiKL(\Lambda; \Sigma) = -\int \Psi \log G^\ast \Lambda G + \tr \Lambda\Sigma.
\end{displaymath}
$\JpsiKL$ is a strictly convex functional over $\LKLG$,
which is an open and convex subset of the Euclidean space $\Rgamma$.
Due to Theorem (\ref{KL_EXISTENCE_THEOREM}),
it does admit a minimum point 
\begin{displaymath}
\LambdaoptKL(\Sigma) = \arg\min_\Lambda \JpsiKL(\Lambda; \Sigma).
\end{displaymath}

\noindent
Let $\delta\Sigma$ be a perturbation of $\Sigma$.
We have
\begin{displaymath}
\begin{split}
\JpsiKL(\Lambda; \Sigma + \delta\Sigma) 
&= -\int \Psi \log G^\ast \Lambda G + \tr \Lambda\Sigma + \tr \Lambda\delta\Sigma\\
&= \JpsiKL(\Lambda; \Sigma) + \left<\delta\Sigma, \Lambda\right>.
\end{split}
\end{displaymath}
It follows from Theorem \ref{CONVEXITYTHEOREM}, where the role of $\delta\Sigma$ is
played by $-p$, that for each $\varepsilon>0$ there exists $\delta>0$ such that
if $||\delta\Sigma||_F<\delta$, then $\JpsiKL(\Lambda; \Sigma + \delta\Sigma)$ 
again admits a minimum point
\begin{equation}
\LambdaoptKL(\Sigma + \delta\Sigma) = \arg\min_\Lambda \JpsiKL(\Lambda; \Sigma + \delta\Sigma)
\end{equation}
and the distance $||\LambdaoptKL(\Sigma + \delta\Sigma) - \LambdaoptKL(\Sigma)||_F$ is less than $\varepsilon$.
The above observation implies well-posedness of the dual problem:
\begin{corollary}
\label{DUALCONTINUITY}
The map
\begin{displaymath}
\Sigma \mapsto \LambdaoptKL(\Sigma)
\end{displaymath}
is continuous from $\SigmaSet$ to $\LKLG$.
\end{corollary}

\noindent
Consider now the primal problem. The variational analysis
yielded the following optimal solution, where the dependence upon
$\Sigma$ has been made explicit:
\begin{displaymath}
\PhioptKL(\Sigma) = \frac{\Psi}{G^\ast\ \LambdaoptKL(\Sigma)\ G}.
\end{displaymath}
We have the following result.
\begin{theorem}
\label{PRIMALCONTINUITY}
The map
\begin{displaymath}
\Sigma \mapsto \PhioptKL(\Sigma)
\end{displaymath}
is a continuous function from $\SigmaSet$ to $L_\infty$.
\end{theorem}

\noindent
\proof
Recall that $\LambdaoptKL(\Sigma)$ is  
the solution of the dual problem where the true asymptotic state variance
is known, and let $\LambdaoptKL(\Sigma + \delta\Sigma)$
be the solution to the dual problem with respect to a perturbed covariance. 
Let $\PhioptKL(\Sigma)$ and $\PhioptKL(\Sigma + \delta\Sigma)$
be the corresponding solutions to the primal problem.
Then
\begin{displaymath}
\begin{split}
||\PhioptKL(\Sigma + \delta\Sigma) - \PhioptKL(\Sigma)||_\infty 
&=    \left\lVert \frac{\Psi}{G^\ast\ \LambdaoptKL(\Sigma + \delta\Sigma)\ G} - \frac{\Psi}{G^\ast\ \LambdaoptKL(\Sigma)\ G}\right\rVert_\infty \\
&\leq ||\Psi||_\infty 
      \left\lVert \frac{1}{G^\ast\ \LambdaoptKL(\Sigma + \delta\Sigma)\ G} - \frac{1}{G^\ast\ \LambdaoptKL(\Sigma)\ G} \right\rVert_\infty. \\
\end{split}
\end{displaymath}
It is easily seen that for each $\eta>0$ we can choose
$\varepsilon >0$ such that if $||\LambdaoptKL(\Sigma + \delta\Sigma) - \LambdaoptKL(\Sigma)||_F < \varepsilon$, then
\begin{eqnarray*}
\nonumber
&&\max_\vartheta |G^\ast\LambdaoptKL(\Sigma + \delta\Sigma) G  - G^\ast \LambdaoptKL(\Sigma) G| 
=\\
\nonumber
&&\ \ \ \ \ \ \ = \max_\vartheta |G^\top(\emjth)(\LambdaoptKL(\Sigma + \delta\Sigma) - \LambdaoptKL(\Sigma)) G(\ejth)| 
<\eta
\end{eqnarray*}

Finally, from the above observation, from Corollary \ref{DUALCONTINUITY},
and from the continuity of the function $\frac{1}{x}$ over $\R^+$, it
follows that for each $\mu>0$, there exists $\delta>0$ such that, for all $||\delta\Sigma||_F < \delta$,
$||\PhioptKL(\Sigma + \delta\Sigma) - \PhioptKL(\Sigma)||_\infty < \mu$.
\qed

\begin{corollary}
\label{HADAMARD}
The problem 
\begin{displaymath}
\arg\min_\Phi D(\Psi||\Phi) \quad \suchthat \quad \int G\Phi G^\ast = \Sigma
\end{displaymath}
is well-posed for $\Sigma\in\SigmaSet$
and for variations $\delta\Sigma$ that belong to $\Rgamma$.
\end{corollary}

\subsection{Well-posedness of Hellinger approximation}

Consider the dual functional (\ref{JPSI_H}):
\begin{displaymath}
\JpsiH(\Lambda; \Sigma) = 
\tr \int (I + G^\ast \Lambda G)^{-1} \Psi + \tr \Lambda\Sigma.
\end{displaymath}
$\JpsiH$ is a strictly convex functional over $\LHG$,
which is an open and convex subset of the Euclidean space $\Rgamma$.
Due to Theorem (\ref{KL_EXISTENCE_THEOREM}),
it admits a minimum point 
\begin{displaymath}
\LambdaoptH(\Sigma) = \arg\min_\Lambda \JpsiH(\Lambda; \Sigma).
\end{displaymath}

\noindent
Let as before $\delta\Sigma$ be a perturbation of $\Sigma$.
Then
\begin{displaymath}
\JpsiH(\Lambda; \Sigma+\delta\Sigma) = \JpsiH(\Lambda; \Sigma) + \left<\delta\Sigma, \Lambda\right>.
\end{displaymath}
Theorem \ref{CONVEXITYTHEOREM} implies the following
\begin{corollary}
\label{DUALCONTINUITY_H}
The map
\begin{displaymath}
\Sigma \mapsto \LambdaoptH(\Sigma)
\end{displaymath}
is continuous from $\SigmaSet$ to $\LHG$.
\end{corollary}

\noindent
The variational analysis yielded the optimal solution
for the primal problem
\begin{equation}
\label{PHIOPTH}
\PhioptH(\Sigma)  = (I + G^\ast\ \LambdaoptH(\Sigma)\ G)^{-1} \Psi (I + G^\ast\ \LambdaoptH(\Sigma)\ G)^{-1},
\end{equation}
and considerations similar to those of theorem 
(\ref{PRIMALCONTINUITY}) lead to the following

\begin{theorem}
\label{PRIMALCONTINUITY_H}
The map
\begin{displaymath}
\Sigma \mapsto \PhioptH(\Sigma)
\end{displaymath}
is continuous from $\SigmaSet$ to $L^{m\times m}_\infty$.
\end{theorem}

\newcommand{\Ql}{Q_\Lambda}
\newcommand{\Qinv}{Q_\Lambda^{-1}}
\newcommand{\Qli}{Q_{\Lambda_i}}
\newcommand{\Qiinv}{Q_{\Lambda_i}^{-1}}

\noindent
To prove  Theorem \ref{PRIMALCONTINUITY_H} we  exploit the following result established
in \cite{RAMPONI_FERRANTE_PAVON_GLOBALLYCONVERGENT} (Lemma 5.2):
\begin{lemma}
\label{LEMMA_RFP} 
Define $\Ql(z) = I + G^\ast(z) \Lambda G(z)$.
Consider a sequence $\Lambda_n\in\LHG$ converging to $\Lambda\in\LHG$. 
Then $Q^{-1}_{\Lambda_n}$ are  well defined and continuous on $\T$
and converge {\em uniformly} to $\Ql^{-1}$ on $\T$.
\end{lemma}

\noindent
\proof (of Theorem \ref{PRIMALCONTINUITY_H}.)
Let $\Ql(z; \Sigma)=I + G^\ast(z)\ \LambdaoptH(\Sigma)\ G(z)$.
Apply Corollary \ref{DUALCONTINUITY_H}
and Lemma \ref{LEMMA_RFP} to establish the continuity
of the map from $\SigmaSet$ to $L^{m\times m}_\infty$
defined by
$\Sigma \mapsto \Ql^{-1}$.
The continuity of $\Sigma \mapsto \PhioptH(\Sigma)$ follows from the continuity
of matrix multiplication.
\qed

\begin{corollary}
\label{HADAMARD_H}
The problem 
\begin{displaymath}
\arg\min_\Phi d_H(\Phi, \Psi) \quad \suchthat \quad \int G\Phi G^\ast = \Sigma
\end{displaymath}
is well-posed, for $\Sigma\in\SigmaSet$
and for variations $\delta\Sigma$ that belong to $\Rgamma$.
\end{corollary}

\section{Consistency}
\label{SECTION_CONSISTENCY}

\noindent
So far we have shown that both the approximation problems
admit an unique solution for all $\Sigma\in\SigmaSet$, and that the
solution is continuous with respect to variations $\delta\Sigma\in\Rgamma$.
The necessity of a restriction to $\Rgamma$ becomes crucial in the
case when we only have an {\em estimate} $\hat{\Sigma}$ of $\Sigma$.

\noindent
In line with the Byrnes-Georgiou-Lindquist theory, and following
an estimation procedure we have sketched in \cite{RAMPONI_FERRANTE_PAVON_GLOBALLYCONVERGENT},
we want to use the above theory to provide an estimate $\hat{\Phi}$ of the
true spectrum of the process $y$.

\noindent
Let $G(z)$ and $\Psi$ be given.
Suppose that we feed $G(z)$ with a {\em finite}
sequence of observations, say $\{y_1, ..., y_N\}$ of the process. Observing the states of the system,
say $\{x_1, ..., x_N\}$, we then compute a Hermitian and positive definite {\em estimate}
$\hat{\Sigma}$ of the asymptotic state covariance, such as
\begin{displaymath}
\hat{\Sigma} = \frac{1}{N} \sum_{k=1}^{N} x_k x_k^\ast.
\end{displaymath}
This is provably consistent, and also unbiased, for we have
supposed from the beginning that $y$ has zero mean.
We seek an estimate $\hat{\Phi}$ of $\Phi$
by solving an approximation problem with respect to $G(z)$, $\Psi$,
and $\hat{\Sigma}$.

\noindent
Since $\hat{\Sigma}$ is not the true variance anymore,
the constraint (\ref{MOMENT_PROBLEM})
may be not feasible.
Hence, in order to find a solution $\hat{\Phi}$,
we need to find a second estimate $\bar{\Sigma}$, close to the first,
such that (\ref{GEORGIOU-CONDITION}) is feasible with the covariance matrix $\bar{\Sigma}$.
%
%
A reasonable way to proceed is to let $\bar{\Sigma}$ be the projection of $\hat{\Sigma}$
onto $\Rgamma$.
Since orthogonal projectors from $\Hermitian(n)$ to a subspace of $\Hermitian(n)$ are
continuous functions, if $\hat{\Sigma}(x_1, ..., x_N)$ is a consistent estimator
of $\Sigma$, then $\bar{\Sigma}$ is also a consistent estimator
of $\Sigma$.

\noindent
The problem that {\em may} come up proceeding in this way is that the projection onto
$\Rgamma$ 
needs
not be positive definite
(that is, it may not belong to $\SigmaSet$), even if $\hat{\Sigma}$ is.
If this is the case, the correct procedure to estimate $\Sigma$ while preserving the structure of
a state covariance compatible with $G(z)$ is 
to find $\bar{\Sigma} \in\SigmaSet$
which is closest to $\hat{\Sigma}$ in a suitable distance.
This is an optimization problem in itself.\\

\noindent
The continuity results of the preceding sections imply two
strong consistency results.
Let $\bar{\Sigma}(x_1, ..., x_N)\in\SigmaSet$ denote a consistent estimator of $\Sigma$.
Let $\PhioptKL(\Sigma)$  be the solution 
to the Kullback-Leibler approximation problem
with respect to the true asymptotic variance
and $\PhioptKL(\bar{\Sigma}(x_1, ..., x_N))$
be the solution of the same problem with respect to the estimate.
%
%
\begin{corollary}
\label{KL_CONSISTENCY}
If   
\begin{equation}
\label{AECONV}
\lim_{N\rightarrow\infty} \bar{\Sigma}(x_1, ..., x_N) = \Sigma \quad{\rm a.s.},
\end{equation}
then 
\begin{displaymath}
\lim_{N\rightarrow\infty} ||\PhioptKL(\bar{\Sigma}(x_1, ..., x_N)) - \PhioptKL(\Sigma)||_\infty = 0 \quad{\rm a.s.}
\end{displaymath}
\end{corollary}
\proof
From the continuity of the map $\Sigma \mapsto \PhioptKL(\Sigma)$ we have that, 
excepting a set of zero probability,
\begin{displaymath}
\lim_{N\rightarrow\infty} \PhioptKL \left( \bar{\Sigma}(x_1(\omega), ..., x_N(\omega)) \right)
= \PhioptKL \left( \lim_{N\rightarrow\infty} \bar{\Sigma}(x_1(\omega), ..., x_N(\omega)) \right)
= \PhioptKL (\Sigma),
\end{displaymath}
where the first limit is taken in $L_\infty(\T)$.
\qed

\noindent
As for the Hellinger multivariable approximation problem,
let $\PhioptH(\Sigma)$ be the solution 
with respect to the true asymptotic variance
and $\PhioptH(\bar{\Sigma}(x_1, ..., x_N))$
be the solution with respect to the estimate. Employing the very same technique used for the proof of Corollary \ref{KL_CONSISTENCY} it is easy to establish the following consistency result for the problem associated to the multivariable Hellinger distance. 
\begin{corollary}
\label{H_CONSISTENCY}
If   
\begin{displaymath}
\lim_{N\rightarrow\infty} \bar{\Sigma}(x_1, ..., x_N) = \Sigma \quad{\rm a.s.},
\end{displaymath}
then 
\begin{displaymath}
\lim_{N\rightarrow\infty} ||\PhioptH(\bar{\Sigma}(x_1, ..., x_N)) - \PhioptH(\Sigma)||_\infty = 0 \quad{\rm a.s.}
\end{displaymath}
\end{corollary}
%

\section{Conclusion}

In this paper, we have considered constrained spectrum approximation problems with respect
to both the Kullback-Leibler pseudo-distance (scalar case) and the  Hellinger distance (multivariable case).
The range of the operator $\Gamma:\Phi\mapsto \int G\Phi G^\ast$ 
is the subspace of the Hermitian matrices that conveyes all the structure that is
needed from a positive-definite matrix in order to be an asymptotic covariance matrix
of the system with tranfer function $G(z)$. As such, it is also a natural subspace to
which the domains of the respective dual problems should be constrained.
We have shown that the condition $\Sigma \in \Rgamma$ is not only necessary for the feasibility
of the moment problem $\{\Phi\ |\ \int G\Phi G^\ast = \Sigma\}$, but also sufficient
for the continuity of the respective solutions with respect to $\Sigma$.
This fact implies well-posedness of both kinds of approximation problems,
and implies the consistency of the respective solutions with respect
to a consistent estimator $\hat{\Sigma}$ of $\Sigma$, as long as it is
restricted to $\Rgamma$. Similar results can be established along the same lines when employing  {\em any} other (pseudo-)distance, as long as the functional form
of the primal optimum depends continuously upon the Lagrange parameter $\Lambda$. 

\bibliographystyle{plain}
\bibliography{fpr}

\end{document}